\documentclass{amsart}
\usepackage[utf8]{inputenc}
\usepackage[T1]{fontenc}

\usepackage{amsmath, amssymb, amsthm, mathtools}
\usepackage{mathrsfs, mathpazo}
\usepackage{geometry}

\usepackage{float}
\usepackage[all]{xy}
\usepackage{graphicx}
\usepackage{color}
\usepackage{tikz}
\usepackage{tikz-cd}
\usepackage{xparse}
\usetikzlibrary{matrix,arrows,decorations.pathmorphing}

\usepackage{cite}
\usepackage{hyperref}
\usepackage{cleveref}

\usepackage{indentfirst}
\usepackage{enumerate}
\usepackage{enumitem}
\usepackage{wrapfig}
\usepackage{adjustbox}
\usepackage{leftidx}
\usepackage{stackrel}
\usepackage{comment}

\newtheorem{theorem}{Theorem}[section]

\newtheorem{lemma}{Lemma}[section]
\newtheorem{remark}{Remark}[section]

\crefname{lemma}{Lemma}{lemma}
\crefname{remark}{Remark}{remark}
\crefname{corollary}{Corollary}{corollary}
\crefname{theorem}{Theorem}{theorem}
\crefname{proposition}{Proposition}{proposition}
\crefname{example}{Example}{example}
\crefname{definition}{Definition}{definition}
\crefname{notation}{Notation}{notation}
\crefname{appendix}{Appendix}{appendix}
\crefname{section}{Section}{section}

\newcommand{\BBB}{\mathcal B}

\newcommand{\HHH}{\mathcal H} %for Hilbert space

\newcommand{\MMM}{\mathcal M}

\newcommand{\MM}{\mathfrak M}

\newcommand{\UUU}{\mathscr U}

\newcommand{\NN}{\mathfrak N}
\newcommand{\JJJ}{\mathcal J}

\newcommand{\PPP}{\mathscr P}
\newcommand{\SSS}{\mathscr S}

\newcommand{\II}{\mathscr I}

\newcommand{\C}{\mathbb C} %for complex number
\newcommand{\R}{\mathbb R}  %for real number
\newcommand{\N}{\mathbb N} % for nature number

\DeclareMathOperator{\re}{Re}

\DeclareMathOperator{\id}{Id}
\DeclarePairedDelimiterX\braket[2]{\langle}{\rangle}{#1 \delimsize\vert #2}

\begin{document}

\title{Surjective $L^p$-isometries of Grassmann spaces}
\author{Wenhua Qian$^1$, Junhao Shen$^2$, Weijuan Shi$^3$, Wenming Wu$^1$, Wei Yuan$^4$}

\address{1. School of Mathematical Sciences, Chongqing Normal University, Chongqing, 601331, China}
\email{whqian86@163.com, wuwm@amss.ac.cn}

\address{2. Department of Mathematics \& Statistics, University of
New Hampshire, Durham, 03824, USA}
\email{Junhao.Shen@unh.edu}

\address{3.School of Mathematics and Information Science,
  Shaanxi Normal University,
  Xi'an , 710062, China}
 \email{shiweijuan@snnu.edu.cn}

\address{4. Institute of Mathematics, Academy of Mathematics and Systems Science \\
Chinese Academy of Sciences, Beijing 100190, China}
\address{4. School of Mathematical Sciences, University of Chinese Academy of Sciences,
Beijing 100049, China}
\email{wyuan@math.ac.cn}

\date{}

\begin{abstract}
Based on the characterization of surjective $L^p$-isometries of unitary groups in finite factors, we describe all surjective $L^p$-isometries between Grassmann spaces of projections with the same trace value in semifinite factors.
\end{abstract}



\thanks{Qian was supported in part by NFS of China (No. 11801150) and Research Foundation of Chongqing Educational Committee(No. KJQN2018000538). Wu was supported in part by NFS of China (No. 11871127, No. 11971463) and Chongqing Science and Technology Commission (No. cstc2019jcyj-msxmX0256). Yuan was supported in part by NFS of China (No. 11871303, No. 11871127, No. 11971463).}
\maketitle

\section{Introduction and statement of the main results}

The research of isometries between substructures of operator algebras has gained much attention recently. One central issue in such a study is investigating if the isometries reflect the algebraic information of the underlying algebra encoded by the substructure, or, more specifically, if the isometries can be described using algebra homomorphisms. The celebrated unitary-antiunitary theorem of Wigner \cite{Wi} is one of the most classical results in this direction. Let $\BBB(\HHH)$ be the algebra of bounded linear operators on the Hilbert space $\HHH$. The Wigner's theorem asserts that every surjective isometry for the operator norm or the $L^2$-norm on the Grassmann space of rank-one projections in $\BBB(\HHH)$ can be extended to a *-automorphism or a *-anti-automorphism of $\BBB(\HHH)$. Due to its importance in the mathematical formulation of quantum mechanics, Wigner's theorem has been studied and generalized by many authors (e.g., \cite{BJM, M1, Pan, Se}).

In recent years, there have been several remarkable developments in the research on isometries of Grassmann spaces. In \cite{Ge1, Ge3}, Geh\'{e}r and $\breve{\mathrm{S}}$emrl described the general form of surjective isometries for the operator norm between Grassmann spaces in $\BBB(\HHH)$. Based on Geh\'{e}r and $\breve{\mathrm{S}}$emrl's work, Mori characterized the surjective isometries for the operator norm between projection lattices of von Neumann algebras using Jordan *-isomorphisms \cite{MM}. As for the isometries with respect to the $L^2$-norm, Geh\'{e}r generalized the early work of Moln\'{a}r \cite{M2} and proved that every (not necessarily surjective) $L^2$-isometry on Grassmann spaces in $\BBB(\HHH)$ is induced by a Jordan *-homomorphisms of $\BBB(\HHH)$ \cite{Ge2}. This result was later generalized to the case of $L^2$-isometry on Grassmann spaces in semifinite factors \cite{QWWY}.

Most recently, Gu, Ma, Shen and Shi studied surjective isometries for the $L^p$-norm on Grassmann spaces in a factor of type $\mathrm{II}$ \cite{Shi}. Let $\MM$ be a semifinite factor equipped with a normal semifinite faithful (n.s.f.) tracial weight $\tau$. For every $A$ in the two-sided ideal $\JJJ$ of operators with finite range projection in $\MM$, its $L^p$-norm is defined by $\|A\|_p := \tau(|A|^p)^{1/p}$ where $0 < p < \infty$. It is well known that $\|\cdot\|_p$ is a norm on $\JJJ$ for $p \in [1, \infty)$. If $p \in (0, 1)$, then $\|\cdot \|_p$ is a quasi-norm on $\JJJ$ such that $\|A_1+A_2\|_p^p \leq \|A_1\|_p^p + \|A_2\|_p^p$ for every $A_1, A_2 \in \JJJ$. We direct the reader to \cite{Tak} for a general reference on the theory of von Neumann algebras. A map $\phi$ between two subsets of $\JJJ$ is called a $L^p$-isometry if
\begin{align*}
    \|B_1 - B_2\|_p = \|\phi(B_1) - \phi(B_2)\|_p
\end{align*}
for every $B_1$, $B_2$ in the domain of $\phi$. The $L^p$-isometry between two subsets in different semifinite factors can be defined similarly. Let $\PPP(\MMM)$ be set of all projections in $\MM$. For $c \in (0, \tau(I))$, we use $\PPP_c(\MM)$ to denote the Grassmann space of all projection in $\MM$ with trace $c$, i.e.,
\begin{align*}
    \PPP_c(\MM) := \{P \in \PPP(\MM): \tau(P) = c\}.
\end{align*}
It is not hard to check that $L^p$-isometries on $\PPP_c(\MM)$ preserve orthogonality in both directions (see \cref{lemma2.6-3-12}). By the structural result for orthogonality preserving maps on $\PPP_c(\BBB(\HHH))$ proved in \cite{Ge1}, we know that every surjective $L^p$-isometry on $\PPP_c(\BBB(\HHH))$ extends to a *-(anti)-automorphism of $\BBB(\HHH)$ if $c \neq \tau(I)/2$. When $\MM$ is a type $\mathrm{II}$ factor and $c \neq \tau(I)/2$, the same result holds and is prove in \cite{Shi}. Although the method used in \cite{Shi} can be modified to give an alternative proof for the case $\MMM = \BBB(\HHH)$, it fails when $\MM$ is a finite factor and $c = \tau(I)/2$. In this paper, we will solve this remaining case and describe the general form of surjective $L^p$-isometries of Grassmann spaces in semifinite factors.

We point out that the surjective $L^p$-isometries on the Grassmann spaces of trace $1/2$ projections in a finite factor are closely related to the surjective $L^p$-isometries of unitary groups. This claim can be backed up by the following observation. Let $\MM = \NN \otimes M_2(\C)$, where $\NN$ is a finite factor. Then the set
\begin{align*}
    \SSS := \left \{ \frac{1}{2} \begin{pmatrix}
        I & U\\
        U^* & I
    \end{pmatrix}: U \in \UUU(\NN) \right\}
\end{align*}
is contained in $\PPP_{1/2}(\MM)$, where $\UUU(\NN)$ is the group of unitaries in $\NN$. If $\phi$ is a surjective $L^p$-isometry on $\PPP_{1/2}(\MM)$ fixing the projection $\begin{psmallmatrix}
    I & 0\\
    0 & 0
\end{psmallmatrix}$, then $\phi$ restricts to a surjective $L^p$-isometry on $\SSS$ and induces a surjective $L^p$-isometry on $\UUU(\NN)$ (see \cref{lem:map_pos}). Therefore, the characterization of surjective $L^p$-isometries on $\UUU(\NN)$ is essential for the study of surjective $L^p$-isometries on $\PPP_{1/2}(\MM)$.

In \cref{sec:II}, we describe the general form of surjective $L^p$-isometries between unitary groups in finite factors. More explicitly, we prove the following theorem, which is one of the main results in this paper.
\begin{theorem}\label{thm:unitary}
    Let $\NN_1$ and $\NN_2$ be two finite factors. If $\psi: \UUU(\NN_1) \to \UUU(\NN_2)$ is a surjective $L^p$-isometry, then there exists a unitary $W \in \NN_2$ and a *-isomorphism or a *-anti-isomorphism $\varphi: \NN_1 \to \NN_2$ such that $\psi$ is of one of the following forms:
\begin{enumerate}
    \item $\psi(U) = W\varphi(U)$ for every $U \in \UUU(\NN_1)$;
    \item $\psi(U)=W\varphi(U^*)$ for every $U \in \UUU(\NN_1)$.
\end{enumerate}
\end{theorem}
\noindent When $\NN_1 = \NN_2 = M_n(\C)$, the proof of \cref{thm:unitary} depends on the methods and ideas in \cite{HM, LM, MS} (see \cref{p-quasi-norm}). For the rest cases, there are two important ingredients in the proof. One is the results on geodesics between unitaries given in \cite{AE} (see \cref{lem:L_2_iso_pre_self_adj}). The other one is the fact that every unitary in $M_{2^n}(\C)$ can be written as the product of self-adjoint unitaries and a complex number of magnitude 1 (see \cref{lem:unitary_by_prod_2_unitary}).

In \cref{sec:III}, base on \cref{thm:unitary}, we describe all surjective $L^p$-isometries on the Grassmann spaces of trace $1/2$ projections of a finite factor. We then reduce the problem of characterizing surjective $L^p$-isometries between Grassmann spaces in different algebras to the characterization of surjective $L^p$-isometries on Grassmann spaces in a factor. Finally, Combining our result with the characterization of surjective $L^p$-isometries on Grassmann spaces given in \cite{Shi}, we prove the following two theorems that describe the general form of surjective $L^p$-isometries between Grassmann spaces of projections with trace $c$.

\begin{theorem}\label{thm:_prop_infty_case}
Let $\phi$ be a surjective $L^p$-isometry from $\PPP_c(\MM_1)$ to $\PPP_c(\MM_2)$, where $c, p \in (0, \infty)$ and $\MM_i$ is a properly infinite semifinite factor equipped with a n.s.f. tracial weight $\tau_i$, $i=1,2$. Then there is a *-isomorphism or a *-anti-isomorphism $\varphi:\MM_1 \to \MM_2$ such that $\tau_1 = \tau_2 \circ \varphi$ and $\phi(P)=\varphi(P)$ for every projection $P\in\PPP_c(\MM_1)$.
\end{theorem}

\begin{theorem}\label{thm:finite_case}
    Let $\phi$ be a surjective $L^p$-isometry from $\PPP_c(\MM_1)$ to $\PPP_c(\MM_2)$, where $c \in (0,1)$, $p \in (0, \infty)$ and $\MM_i$ is a finite factor equipped with a tracial state $\tau_i$, $i=1,2$. If $c \neq 1/2$, there is a *-isomorphism or a *-anti-isomorphism $\varphi:\MM_1 \to \MM_2$ such that $\phi$ is of the following form
    \begin{align}\label{equ:finite_case_1}
      \phi(P)=\varphi(P), \quad \forall P\in\PPP_c(\MM_1).
   \end{align}
In the case when $c = 1/2$, we have either (\ref{equ:finite_case_1}), or the following additional possibility occurs:
\begin{align}
    \phi(P) = I- \varphi(P), \quad \forall Q \in \PPP_{1/2}(\MM_1).
\end{align}
\end{theorem}

\section{surjective $L^p$-isometries of unitaries in finite factors}\label{sec:II}

This section is devoted to proving \cref{thm:unitary}. Let $\NN$ be a finite factor with the tracial state $\tau$. We use $\UUU(\NN)$ to denote the set of unitaries in $\NN$. In \cite{LM}, Moln\'{a}r showed that every isometry on $\UUU(M_n(\C))$ with respect to a unitarily invariant norm is induced by a *-automorphism or a *-anti-automorphism on $M_n(\C)$. For $p \in [1, \infty)$, $\|\cdot\|_p$ is a unitarily invariant norm on $\UUU(M_n(\C))$. We thus know that every $L^p$-isometry $\psi$ on $\UUU(M_n(\C))$ satisfying $\psi(I) = I$ can be extended to a *-automorphism or a *-anti-automorphism if $p \geq 1$. We now show that the same result holds for $L^p$-isometries on $\UUU(M_n(\C))$ when $0 < p < 1$.

\begin{lemma}\label{p-quasi-norm}
Let $\psi$ be a surjective $L^p$-isometry on $\mathscr{U}(M_n(\mathbb{C}))$. Then there exists a unitary $W \in \mathscr{U}(M_n(\mathbb{C}))$ and a *-automorphism or a *-anti-automorphism $\varphi$ on $M_n(\C)$ such that $\psi$ is of one of the following forms:
\begin{enumerate}
    \item $\psi(U)= W\varphi(U)$, $U\in\mathscr{U}(M_n(\mathbb{C}))$;
    \item $\psi(U)= W\varphi(U^*)$, $U\in\mathscr{U}(M_n(\mathbb{C}))$.
\end{enumerate}
\end{lemma}

\begin{proof}
Without loss of generality, we may assume that $\psi(I) = I$. Otherwise, we can consider the map $\psi(I)^* \psi(\cdot)$ instead. The lemma will be proved if we can show that there exists a *-automorphism or a *-anti-automorphism $\varphi$ such that $\psi(U) = \varphi(U)$ for every unitary $U$ or $\psi(U) = \varphi(U^*)$ for every unitary $U$. By Theorem 3 in \cite{LM}, we only need to prove this result for the case $0 < p < 1$.

Since $\|A\|=\||A|\|$, we have
\begin{align*}
    \|A\|_p \leq \|A\| \leq \sqrt[p]{n}\|A\|_p, \quad \forall A \in M_n(\C).
\end{align*}
In particular, $\psi$ is continuous with respect to the operator norm. By Proposition 10 in \cite{LM}, we only need to show that $\psi$ preserves the inverted Jordan triple product, i.e.,
    \begin{align}\label{equ:p-quasi-norm_eq_1}
        \psi(VW^*V) = \psi(V) \psi(W)^* \psi(V), \quad \forall V, W \in \UUU(M_n(\C)).
    \end{align}

    Let $c:=\frac{1}{\sqrt[p]{n}}$ and $V, W \in \UUU(M_n(\C))$ such that $\|V - W\| \leq c^2/2$. We claim that
    \begin{align}
        \sqrt{3} \|VX^* - I\|_p \leq \|VX^*VX^* - I\|_p
    \end{align}
    for every unitary $X$ satisfying
    \begin{align*}
        \|X-W\|_p=\|VW^*V-X\|_p=\|V-W\|_p.
    \end{align*}
Let $U = VX^*$. Note that
\begin{align*}
\|U-I\|_p&=\|V-X\|_p \leq 2 \max\{\|V-W\|,\|W-X\|\} \leq \frac{2}{c}\|V-W||_p\leq c.
\end{align*}
Therefore $\|U - I\| \leq 1$. This implies that the spectrum of $U$ is contained in the set $\{\lambda:|\lambda|=1, |\lambda-1|\leq 1\}$. Then it is not hard to see that $U + I$ is invertible and $\|(U+ I)^{-1}\| \leq 1/ \sqrt{3}$. Consequently, we have
\begin{align*}
    &\sqrt{3}\|U-I\|_p = \sqrt{3}\|(U+I)^{-1}(U^2-I)\|_p  \leq \sqrt{3} \|(U+I)^{-1}\|\|U^2-I\|_p \leq \|U^2-I\|_p,
\end{align*}
and the claim is proved.

By Corollary 3.9 in \cite{HHMM}, we know that the equation (\ref{equ:p-quasi-norm_eq_1}) holds for every unitaries $V, W$ satisfying $\|V - W\| \leq c^2/2$. Proceed as in the proof of Theorem 8 in \cite{HM}, we can show that the equation (\ref{equ:p-quasi-norm_eq_1}) holds for every pair of unitaries. Hence the lemma is proved.
\end{proof}

\begin{remark} As $\mathscr{U}(M_n(\mathbb{C}))$ is a compact connected manifold, then according to the invariance domain theorem the $L^p$-isometry on $\mathscr{U}(M_n(\mathbb{C}))$ must be surjective. We give the surjective assumption in the above Lemma for the consistency.
\end{remark}

Let $\NN_i$ be a finite factor with the traical state $\tau_i$, $i=1,2$. We next show that $\NN_1$ and $\NN_2$ are either *-isomorphic or *-anti-isomorphic if there exists a $L^p$-isometry from $\UUU(\NN_1)$ onto $\UUU(\NN_2)$. We prepare the proof by establishing three lemmas.

\begin{lemma} \label{lemma2.4-3-12}
Let $U\in\mathscr U(\NN)$, where $\NN$ is a finite factor. Assume that $p \neq 2$. Then $\|I-U\|_p^p+\|I+U\|_p^p=2^p$ if and only if $U=U^*$.
\end{lemma}

\begin{proof}
If $\|I-U\|_p^p+\|I+U\|_p^p=2^p$, then
\begin{align*}
 2(\|I-U\|_p^p+\|I+U\|_p^p) = \|(I-U)+(I+U)\|_p^p+\|(I-U)-(I+U)\|_p^p.
\end{align*}
By Corollary 2.1 in \cite{Lam}, we have $(I+U)^*(I-U)=0$. Thus $U=U^*$.

Conversely, if $U=U^*$, then
\begin{align*}
 \|I-U\|_p^p+\|I+U\|_p^p= 2^p\tau (P)+2^p\tau (Q)=2^p,
\end{align*}
where $P=(I+U)/2$ and $Q=(I-U)/2$ are orthogonal projections.
\end{proof}

\begin{lemma}\label{lem:L_2_iso_pre_self_adj}
Let $\psi: \UUU(\NN_1) \to \UUU(\NN_2)$ be a surjective $L^p$-isometry such that $\psi(I) = I$, where $\NN_1$ and $\NN_2$ are two finite factors. For every $U \in \UUU(\NN_1)$, we have $U^* = U$ if and only if $\psi(U)^* = \psi(U)$.
\end{lemma}

\begin{proof}
    If $p \neq 2$, the lemma is proved by \cref{lemma2.4-3-12}. Assume that $p = 2$. Since $\psi$ is a $L^2$-isometry, we have
\begin{align*}
    \re(\tau_1(U^*V))=\re(\tau_2(\psi(U)^*\psi(V))), \quad \forall \quad U,V\in\mathcal{U}(\NN).
\end{align*}
Thus it is easy to check that $\|\sum_{i=1}^n\lambda_iU_i\|_2=\|\sum_{i=1}^n\lambda_i\psi(U_i)\|_2$, where $\lambda_i \in \R$, $U_i \in \UUU(\NN_1)$. Recall that for every $A \in \NN_1$, there exists two unitaries $W_1$, $W_2$ in $\NN_1$ such that $A = \frac{\|A\|}{2}(W_1 + W_2)$. Then $\psi$ has a unique real-linear extension $\Psi:\NN_1 \to \NN_2$, i.e.,
\begin{align*}
    \Psi(\sum_{i=1}^n\lambda_i U_i) :=\sum_{i=1}^{n}\lambda_i\psi(U_i).
\end{align*}
It is clear that $\Psi$ is a surjective $L^2$-isometry.

Let $U= 2P-I$, where $P$ is a projection in $\NN_1$. Then $\psi(U)^* = \psi(U)$ if and only if $\Psi(P)^* = \Psi(P)$. By Lemma 3.1 in \cite{AE}, the curve $t \to e^{itP}$, $t \in [0, 1]$, which is $C^\infty$ as a curve in the Hilbert space $L^2(\NN_1, \tau_1)$, has minimal length among piecewise C$^1$ curves of unitaries jointing $I$ and $e^{iP}$ in the $L^2$ metric. Since $\Psi$ is a real-linear $L^2$-isometry, $f(t):=\psi(e^{itP})$ is a C$^\infty$ curve in $L^2(\NN_2, \tau_2)$ with
    \begin{align*}
        \frac{df}{dt}(t) = \Psi(iPe^{itP}), \quad \frac{d^2f}{dt^2}(t) = -\Psi(Pe^{itP}).
    \end{align*}
Since $\Psi$ is bijective, the curve $f$ has minimal length among piecewise C$^1$ curves of unitaries jointing $I$ and $\psi(e^{iP})$. By Theorem 2.4 and Proposition 2.2 in \cite{AE}, there exists a self-adjoint operator $H \in \NN_2$ such that $f(t) = e^{itH}$. Therefore $\Psi(P) = -\frac{d^2f}{dt^2}(0) = H^2$. In particular, $\Psi(P)$ is self-adjoint.
\end{proof}

\begin{lemma}\label{lemma2.6-3-12}
Let $\MM$ be a semifinite factor with a n.s.f. tracial weight $\tau$ and $P_1, P_2$ be two finite projections in $\MM$, i.e., $\tau(P_i) < \infty$. Then $P_1\perp P_2$ if and only if
\begin{align*}
\|P_1-P_2\|_p^p=\tau(P_1)+\tau(P_2).
\end{align*}
\end{lemma}

\cref{lemma2.6-3-12} can be proved by the method analogous to that used in the proof of Lemma 2.4 in \cite{Shi}, and its proof is omitted.

\begin{theorem}\label{prop2.7-3.14}
Let $\psi: \UUU(\NN_1) \to \UUU(\NN_2)$ be a surjective $L^p$-isometry such that $\psi(I) = I$, where $\NN_1$ and $\NN_2$ are two finite factors. Then there exists a *-isomorphism or *-anti-isomorphism $\varphi: \NN_1 \to \NN_2$ such that $\psi(U) = \varphi(U)$ for every self-adjoint unitary $U$.
\end{theorem}

\begin{proof}
     Let
\begin{align*}
 \phi(P) :=\frac {\psi(2P-I)+I}{2}, \qquad \forall  P\in \mathscr P(\NN_1).
\end{align*}
By \cref{lem:L_2_iso_pre_self_adj}, $\phi$ is a $L^p$-isometry $\phi$ on $\mathscr P(\NN_1)$. Since $\psi$ is surjective, $\phi$ is also surjective. Note that $\phi(I)=I$, we have
\begin{align*}
 1-\tau_1(P)=\|I-P\|_p^p=\|I-\phi(P)\|_p^p =1-\tau_2(\phi(P)),
\end{align*}
    where $\tau_i$ is the tracial state on $\NN_i$. Therefore $\tau_2(\phi(P))=\tau_1(P)$. By \cref{lemma2.6-3-12}, $\phi$ preserves orthogonality in both directions.

     If $\NN_1\ncong M_2(\C)$, then there exists a *-isomorphism or a *-anti-isomorphism $\varphi: \NN_1 \to \NN_2$ such that $\phi(P) = \varphi(P)$ for every $P\in \mathscr P(\NN)$ by Corollary in \cite{Dye} (see also \cite{Kad}). Thus $\psi(U) = \varphi(U)$ for every self-adjoint unitary.

    If $\NN_1 \cong M_2(\C)$, then we have $\NN_2\cong M_2(\C)$. And \cref{p-quasi-norm} implies the result.
\end{proof}

By \cref{p-quasi-norm} and \cref{prop2.7-3.14}, we only need to prove \cref{thm:unitary} under the assumption that $\NN_1 = \NN_2 =\NN$, where $\NN$ is a $\mathrm{II}_1$ factor, and $\psi$ is a surjective $L^p$-isometry on $\UUU(\NN)$. From now on, we use $\NN$ to denote a $\mathrm{II}_1$ factor and $\psi$ to denote a surjective $L^p$-isometry on $\UUU(\NN)$. The proof of \cref{thm:unitary} is divided into two cases: $p = 2$ and $p \neq 2$.

\subsection{The $p = 2$ case}

\begin{proof}[\bf{Proof of \cref{thm:unitary}}($p = 2$)]
By \cref{prop2.7-3.14} and \cref{p-quasi-norm}, we only need to consider the case that $\NN_1 = \NN_2 =\NN$, where $\NN$ is a II$_1$ factor. Without loss of generality, we may assume that $\psi(I) = I$. Otherwise, we can consider the map $\psi(I)^* \psi(\cdot)$ instead. By \cref{prop2.7-3.14}, there exists a *-automorphism or a *-anti-automorphism $\varphi$ on $\NN$ such that $\psi(U) = \varphi(U)$ for every self-adjoint unitary $U$. By considering the map $\varphi^{-1}\circ\psi$ instead, we may further assume that $\psi(U) = U$ for every self-adjoint unitary $U$.

By the proof of \cref{lem:L_2_iso_pre_self_adj}, we know $\psi$ has a unique real-linear extension $\Psi$ which is a surjective $L^2$-isometry on $\NN$. Moreover, $\Psi(P)=P$ for every projection $P$ and
    \begin{align*}
        \re \tau(A^*B) = \re \tau(\Psi(A)^*\Psi(B)),
    \end{align*}
for every pair of elements $A, B \in \NN$, where $\tau$ is the tracial state on $\NN$. It is clear that the conclusion will be proved if we can show that $\Psi(A) = A$ for every $A \in \NN$ or $\Psi(A) = A^*$ for every $A \in \NN$.

Note that $\Psi((U + iI)/\sqrt{2}) = (U + \psi(iI))/\sqrt{2}$ is a unitary for every self-adjoint unitary $U$. Thus $\psi(iI)^* = -\psi(iI)$ and $\psi(iI)U = U \psi(iI)$. Since $\mathfrak{N}$ is a factor, $\psi(iI) = iI$ or $-iI$. We claim that $\Psi(iP)= \psi(iI)P$ for every projection $P$ in $\NN$. Since
\begin{align*}
    (I-P) + \Psi(iP) = \psi((I-P)+iP) \mbox{ and }  (I-P) - \Psi(iP) = \psi((I-P)-iP)
\end{align*}
are unitaries, we have $\Psi(iP)^*\Psi(iP) = \Psi(iP)\Psi(iP)^* = P$. Note that
\begin{align*}
    \tau(P) = \re \tau(\psi(iI)^* \psi((I-P) + iP)) = \re \tau(\psi(iI)^* \Psi(iP)).
\end{align*}
Therefore $\Psi(iP) = \psi(iI)P$. Since every operator can be approximated by linear combinations of projections, we have $\Psi(A) = A$ for every $A\in \NN$ or $\Psi(A) = A^*$ for every $A \in \NN$.
\end{proof}

\subsection{The $p \neq 2$ case}
In the rest of this subsection, we assume that $p \neq 2$.

\begin{lemma}\label{lem:scalar_map}
    Let $U \in \UUU(\NN)$ and $\theta \in (-\pi, \pi]$. Then
    \begin{align}\label{equ:scalar_map_1}
        \|U - 2P+I\|_p^p = \tau(P)(2-2\cos \theta)^{p/2} +  \tau(I-P)(2+2\cos \theta)^{p/2}
    \end{align}
    for every projection $P \in \NN$ if and only if $U$ equals either $e^{i\theta}I$ or $e^{-i\theta}I$.
\end{lemma}

\begin{proof}
We only prove the sufficiency. Note that $-I$ is the only unitary satisfying  $\|U - I\|_p = 2^p$ and $I$ is the only unitary satisfying  $\|U + I\|_p = 2^p$. Then the the statement of the lemma is true if $\theta = 0$ or $\pi$.

Let $\theta \in (-\pi, 0) \cup (0, \pi)$. We need to show that the equation $(\ref{equ:scalar_map_1})$ implies that $U = e^{i\theta}I$ or $e^{-i\theta}I$. Without loss of generality, we could assume that $0 < \theta < \pi$. Let $\chi$ be the characteristic function of the set $\{e^{it}: t \in [\theta, 2\pi-\theta]\}$. Note that $1-\cos t_1 > 1-\cos \theta$ for every $t_1 \in (\theta, 2\pi-\theta)$ and $1+ \cos t_2 > 1+ \cos \theta$ for every $t_2 \in (-\theta, \theta)$. The equation
\begin{align*}
    \|U - 2\chi(U)+I\|_p^p = \tau(\chi(U))(2-2\cos \theta)^{p/2} +  \tau(I-\chi(U))(2+2\cos \theta)^{p/2}
\end{align*}
implies that $U = e^{i\theta}E + e^{-i\theta}(I-E)$, where $E$ is a projection in $\NN$. We now show that $E = 0$ or $I$.

Assume that $E \neq 0$ and $E \neq I$. We only consider the case that $E \preceq I-E$. And the proof for the case $I-E \preceq E$ is similar, so is omitted. Let $W$ be a partial isometry such that $W^*W = E$ and $WW^* \leq I-E$. Since $V = (W+ W^*) - (I-E- WW^*)$ is a self-adjoint unitary, we have
    \begin{align*}
        \tau(E)(2-2\cos \theta)^{\frac{p}{2}} +  \tau(I-E)(2+2\cos \theta)^{\frac{p}{2}} &=
        \|U - V\|_p^p \\
        &= \tau(I -2E)(2+2\cos\theta)^{\frac{p}{2}} + \tau(E)2^p.
    \end{align*}
Therefore $(2-2\cos \theta)^{\frac{p}{2}} + (2+2\cos \theta)^{\frac{p}{2}} = 2^p$. Recall that $p \neq 2$ and $|\cos \theta | < 1$, we have
    \begin{align*}
        \left [(\frac{1-\cos \theta}{2})^{\frac{p}{2}} + (\frac{1+\cos \theta}{2})^{\frac{p}{2}} \right ] \neq 1.
    \end{align*}
This leads to a contradiction.
\end{proof}

\begin{lemma}\label{lem:desc_scalar}
    Let $\psi$ be a surjective $L^p$-isometry on $\UUU(\NN)$ such that $\psi(I) = I$. Then we have either $\psi(e^{i\theta}I) = e^{i\theta}I$ for every $\theta \in \R$ or $\psi(e^{i\theta}I) = e^{-i\theta}I$ for every $\theta \in \R$.
\end{lemma}

\begin{proof}
By \cref{prop2.7-3.14}, we may assume that $\psi(U) = U$ for every self-adjoint unitary. Then \cref{lem:scalar_map} and \cref{p-quasi-norm} implies the result.
\end{proof}

\begin{lemma}\label{lem:phase_inv}
    Let $\psi$ be a surjective $L^p$-isometry on $\UUU(\NN)$. If $\psi(U) = U$ for every self-adjoint unitary and $\psi(e^{i\theta}I) = e^{i\theta}I$ for every $\theta \in \R$, then $\psi(e^{i\theta}U) = e^{i\theta}U$ for every self-adjoint unitary $U$ and $\theta$.
\end{lemma}

\begin{proof}
    Let $P$ be a projection in $\NN$. Note that $\psi'(\cdot) : = (2P-I)\psi((2P-I)(\cdot))$ is a surjective $L^p$-isometry on $\UUU(\NN)$ satisfying $\psi'(I) = I$. By \cref{lem:desc_scalar}, we have either $\psi(e^{i\theta}(2P-I)) = e^{i\theta}(2P-I)$ for every $\theta \in \R$ or $\psi(e^{i\theta}(2P-I)) = e^{-i\theta}(2P-I)$ for every $\theta \in \R$. If we can show that $\psi(i(2P-I)) = i(2P-I)$, the lemma is proved.

Assume that $\psi(i(2P-I)) = -i(2P-I)$. Let $Q$ be a projection in $\NN$. Note that
\begin{align*}
   \begin{cases}
       \|P-Q\|_p + \|P+Q -I\|_p \geq \|2P -I\|_p = 1, \quad & p \geq 1\\
       \|P-Q\|_p^p + \|P+Q-I\|_p^p \geq \|2P-I\|_p^p = 1, \quad &0 < p < 1
   \end{cases}
\end{align*}
    (see Chapter 5 in \cite{XBC}). Thus we have $\psi(i(2Q-I)) = -i(2Q-I)$ if $\|P-Q\|_p$ is small enough. Since the set of projections in $\NN$ is connected with respect to the topology induced by the by the quasi-norm $\|\cdot\|_p$, we have $\psi(i(2E - I)) = -i(2E-I)$ for every projection $E$ in $\NN$. Recall that $\psi(iI) = iI$. We have a contradiction. Therefore $\psi(i(2P-I)) = i(2P-I)$.
\end{proof}

\begin{lemma}\label{product_two_unitaries_1}
    Let $\psi$ be a surjective $L^p$-isometry on $\mathscr{U}(\NN)$ such that $\psi(e^{i\theta}U)=e^{i\theta}U$ for every self-adjoint unitary $U\in \NN$ and $\theta \in \R$. Then we have
   \begin{align*}
        \psi(e^{i\theta}UV)= e^{i\theta}UV
   \end{align*}
for every self-adjoint unitary $U, V \in \NN$.
\end{lemma}

\begin{proof}
    Let $P$ be a trace $1/2$ projection. Then $(2P-I)\psi((2P-I)(\cdot))$ is a surjective $L^p$-isometry on $\UUU(\NN)$. Note that $(2P-I)\psi(e^{i\theta}(2P-I)) = e^{i\theta}I$. By \cref{prop2.7-3.14} and \cref{lem:phase_inv}, there exists a *-automorphism or a *-anti-automorphism $\varphi$ on $\NN$ such that $(2P-I) \psi(e^{i\theta}(2P-I)V) = e^{i\theta}\varphi(V)$ for every self-adjoint unitary $V$ and $\theta \in \R$. We claim that $\varphi$ is the identity *-isomorphism. Indeed, note that
$\varphi(P) = P$ and
\begin{align*}
    \varphi(W) &=\frac{1}{2}(\varphi(W + W^*) + \varphi(W - W^*))\\
    &= \frac{1}{2}(2P-I)(\psi(W-W^*) + \psi(W+W^*)) = W,
\end{align*}
    for every partial isometry $W$ satisfying $W^*W= I-P$ and $WW^* = P$. We have $\varphi(PW)= \varphi(P)\varphi(W)$. Thus $\varphi$ is a *-isomorphism. Since $\NN = \{W, W^*: W\in \NN, W^*W= I-P, WW^* = P\}''$, $\varphi$ is the identity *-isomorphism. Therefore we conclude that
    \begin{align*}
     \psi(e^{i\theta}(2P-I)V) = e^{i\theta}(2P-I)V
    \end{align*}
for every trace $1/2$ projection $P$ and self-adjoint unitary $V$.

    Since $\psi'(\cdot): =e^{-i\theta}\psi(e^{i\theta}(\cdot)V)V$ is a surjective $L^p$-isometry such that $\psi'(I) = I$, there exists a *-automorphism or a *-anti-automorphism $\varphi'$ on $\NN$ such that $\psi'(U) = \varphi'(U)$ for every self-adjoint unitary $U \in \NN$ by \cref{prop2.7-3.14}. Based on the above discussion, we have
    \begin{align*}
        \varphi'(2P-I) = \psi'(2P-I) = 2P-I
    \end{align*}
for every trace $1/2$ projection $P$ in $\NN$. Recall that the linear span of the set of trace $1/2$ projections is dense in $\NN$ (see Lemma 3.7 in \cite{QWWY}). Therefore, $\varphi'$ is the identity *-isomorphism, and the lemma is proved.
\end{proof}

\begin{lemma}\label{product_sa_unitaries}
    Let $\psi$ be a surjective $L^p$-isometry on $\mathscr{U}(\NN)$ such that $\psi(U)=U$ for every self-adjoint unitary $U \in \NN$ and $\psi(e^{i\theta}I) = e^{i\theta}I$ for every $\theta \in \R$. Then we have
\begin{align*}
    \psi(e^{i\theta}U_1 \cdots U_n)=e^{i\theta}U_1 \cdots U_n
\end{align*}
for self-adjoint unitaries $U_1, \ldots, U_n \in \NN$ and $\theta \in \R$, $n \geq 1$.
\end{lemma}

\begin{proof}
    We prove the lemma by induction. If $n=2$, the lemma holds by \cref{lem:phase_inv} and \cref{product_two_unitaries_1}. Assume that the lemma is true for $k \geq 2$. Let $V = U_1 \cdots U_{k-1}$, where $U_i$ is a self-adjoint unitary in $\NN$, $i=1,\ldots,k-1$. Note that $\psi'(\cdot):=V^*\psi(V(\cdot))$ is a surjective $L^p$-isometry such that $\psi'(e^{i\theta}U_k) = e^{i\theta}U_k$ for every self-adjoint unitary $U_k$. By \cref{product_two_unitaries_1}, we have
    \begin{align*}
        (U_1 \cdots U_{k-1})^* \psi(e^{i\theta}U_1 \cdots U_{k+1})=e^{i\theta}U_{k}U_{k+1}.
    \end{align*}
And the lemma is proved.
\end{proof}

The following lemma is well-known to experts, and we include the proof for the sake of completeness.

\begin{lemma}\label{lem:unitary_by_prod_2_unitary}
    Every unitary $V$ in $M_{2^n}(\C) (n\geq 1)$ can be written as
    \begin{align*}
        V = e^{i \theta} U_1 \cdots U_l
    \end{align*}
    where $U_i$ are self-adjoint unitaries in $M_{2^n}(\C)$.
\end{lemma}

\begin{proof}
    Let $E_{ij}$, $i,j = 1, \ldots, 2^n$, be the canonical matrix units of $M_{2^n}(\C)$. We only need to show that there exist self-adjoint unitaries $U_1, \ldots, U_{2n}$ and $\theta \in \R$ such that
    \begin{align*}
        e^{i\beta} E_{11} + (I-E_{11}) = e^{i \theta} U_1 \cdots U_{2n},
    \end{align*}
    for every $\beta \in \R$.
We prove this fact by induction on $n$.

For $n = 1$, we have
\begin{align}\label{equ:unitary_by_prod_2_unitary_1}
    e^{i\beta} E_{11} + E_{22} = e^{\frac{i\beta}{2}}(iE_{12} - iE_{21})(ie^{\frac{-i\beta}{2}}E_{12} -ie^{\frac{i\beta}{2}}E_{21}).
\end{align}
Assume that the fact is true for $k \geq 1$. Let $m = 2^k$. Then there exist self-adjoint unitaries $U_1, \ldots, U_{2k} \in M_{2^{k+1}}(\C)$ and $\theta \in \R$ such that
\begin{align*}
    e^{i\beta} E_{11} + (E_{22} + \cdots + E_{mm}) + e^{i\theta}(I- E_{11} - \cdots - E_{mm}) = e^{i\theta} U_1 \cdots U_{2k}.
\end{align*}
By equation (\ref{equ:unitary_by_prod_2_unitary_1}), there exist two self-adjoint unitaries $W_1, W_2 \in M_{2^{k+1}}(\C)$ such that
    \begin{align*}
        (E_{11} + E_{22} + \cdots + E_{mm}) + e^{-i\theta}(I- (E_{11} + E_{22} + \cdots + E_{mm})) = e^{\frac{-i\theta}{2}}W_1W_2.
    \end{align*}
Therefore $e^{i\beta} E_{11} + (I-E_{11}) = e^{\frac{i\theta}{2}} U_1 \cdots U_{2k} W_1W_2$.
\end{proof}

\begin{proof}[\bf{Proof of \cref{thm:unitary}} ($p \neq 2$)]
By \cref{prop2.7-3.14} and \cref{p-quasi-norm}, we only need to consider the case that $\NN_1 = \NN_2 =\NN$, where $\NN$ is a II$_1$ factor. Without loss of generality, we may assume that $\psi(I) = I$. By \cref{prop2.7-3.14}, there exists a *-automorphism or a *-anti-automorphism $\varphi$ on $\NN$ such that $\psi(U) = \varphi(U)$ for every self-adjoint unitary $U$. Let
    \begin{align*}
        \psi'(V) := \begin{cases}
            \varphi^{-1} \circ \psi(V), \quad & \mbox{if $\psi(e^{i \theta}I) = e^{i \theta}I$}\\
            \varphi^{-1} \circ \psi(V^*), \quad & \mbox{if $\psi(e^{i \theta}I) = e^{-i \theta}I$}\\
        \end{cases}
    \end{align*}
    By \cref{lem:phase_inv}, $\psi'$ is a surjective $L^p$-isometry such that $\psi'(e^{i \theta}U) = e^{i \theta}U$ for every self-adjoint unitary $U$. We now show that $\psi'$ is the identity map on $\UUU(\NN)$.

    Let $V \in \UUU(\NN)$. For every $k \in \N$, there exist an integer $n$ and a subfactor $\NN_k$ of $\NN$ such that $\NN_k \simeq M_{2^{n}}(\C)$ and a unitary $V_k \in \NN_k$ satisfying $\|V - V_k\|_p < 1/k$. By \cref{product_sa_unitaries} and \cref{lem:unitary_by_prod_2_unitary}, $\psi'(V_k) = V_k$. Since
   \begin{align*}
      \lim_k \|\psi'(V) - V_k\|_p = \lim_k \|V - V_k\|_p  = 0,
   \end{align*}
we have $\psi'(V) = \lim_{k} V_k = V$.
\end{proof}

\section{Surjective $L^p$-isometries of the Grassmann spaces}\label{sec:III}

Let $\MM$ be a semifinite factor equipped with a n.s.f. tracial weight $\tau$. If $\MM$ is a finite factor, we require that $\tau$ is the tracial state on $\MM$, i.e., $\tau(I) = 1$. In the following, we always assume that $\PPP_c(\MM) \neq \emptyset$. Let $\phi$ be a surjective $L^p$-isometry on $\PPP_c(\MM)$. Then $\phi$ preserves the orthogonality in both direction by \cref{lemma2.6-3-12}. By Theorem 1.2 in \cite{Ge1} and Theorem 4.9 in \cite{Shi}, $\phi$ can be extended to a *-automorphism or a *-anti-automorphism on $\MM$ if $c \neq \tau(I)/2$. If $\MM$ is a finite factor and $c = 1/2$, it is clear that $\varphi(\cdot)$ and $I-\varphi(\cdot)$ are both surjective $L^p$-isometries on $\PPP_{1/2}(\MM)$, where $\varphi$ is a *-automorphism or a *-anti-automorphism on $\MM$. We now show that every surjective $L^p$-isometry on $\PPP_{1/2}(\MM)$ arises in this manner.

\begin{lemma} \label{lemma3.3-3-14}
Let $\MM$ be a finite factor equipped with a tracial state $\tau$ and $P, Q\in\PPP_{\frac{1}{2}}(\mathfrak{M})$. If $p \neq 2$, then
\begin{align*}
   \|P -Q \|_{p}=\|(I-P) - Q \|_{p}=\frac {1}{\sqrt 2}
\end{align*}
if and only if there exists a partial isometry $U\in\mathfrak{M}$ such that $U^*U=P$, $UU^*=I-P$ and $Q=(I+U+U^*)/2$.
\end{lemma}

\begin{proof}
    Without loss of generality, we may assume that $\MM = \NN \otimes M_2(\mathbb{C})$ and
\begin{align*}
    P = I \otimes E_{11},\quad Q = H \otimes E_{11} + \sqrt{H(I-H)} \otimes (E_{12} + E_{21}) + (I-H) \otimes E_{22},
\end{align*}
    where $H$ is a positive contraction in the finite factor $\NN$ and $\{E_{ij}\}$ is the canonical matrix units of $M_2(\C)$. Note that
\begin{align*}
    \|I-H\|_{p/2} = \|P-Q\|_p^2 = \|(I-P)-Q\|_p^2 = \|H\|_{p/2} =1/2,
\end{align*}
where $\|\cdot\|_{p/2}$ is the $L^{p/2}$-norm on $\NN$ induced by the tracical state on $\NN$. Therefore
\begin{align*}
\|(I-H)+H\|_{p/2}=\|I-H\|_{p/2}+\|H\|_{p/2}.
\end{align*}
    Note that $p/2 \neq 1$. By the condition for Minkowski's inequality to be an equality (see 6.13 in \cite{HLP}), we know that there exists a positive number $\lambda$ such that $H=\lambda (I-H)$. Since $\|I-H\|_{p/2} = \|H\|_{p/2}$, $\lambda = 1$. And the lemma is proved.
\end{proof}

\begin{lemma}\label{lem:map_pos}
    Let $\MM_i = \NN_i \otimes M_2(\C)$, where $\NN_i$ is a finite factor, $i=1,2$. If $\phi$ is a surjective $L^p$-isometry from $\PPP_{1/2}(\MM_1)$ to $\PPP_{1/2}(\MM_2)$ satisfying $\phi(I \otimes E_{11}) = I \otimes E_{11}$, then there exists a surjective $L^p$-isometry $\psi: \UUU(\NN_1) \to \UUU(\NN_2)$ such that
    \begin{equation}\label{equ:map_pos}
 \begin{aligned}
    &\phi(\frac{1}{2}(I \otimes (E_{11} + E_{22}) + V \otimes E_{12} + V^* \otimes E_{21}))\\
    =&\frac{1}{2}(I \otimes (E_{11} + E_{22}) + \psi(V) \otimes E_{12} + \psi(V)^* \otimes E_{21})
\end{aligned}
    \end{equation}
    for every unitary $V \in \NN_1$, where $\{E_{ij}\}_{ij}$ is the canonical matrix units of $M_2(\C)$.
\end{lemma}

\begin{proof}
We only need to show that for every unitary $V \in \NN_1$, there exists a unitary $\psi(V) \in \NN_2$ such that the equation (\ref{equ:map_pos}) holds. By \cref{lemma2.6-3-12}, $\phi(I \otimes E_{22}) = I \otimes E_{22}$. If $p \neq 2$, then \cref{lemma3.3-3-14} implies the result.

Assume that $p = 2$. Proceed as in the proof of Lemma 1 in \cite{M2} (or see the proof of Lemma 2.2 in \cite{QWWY}), we can show that $\phi$ has a unique linear extension $\Phi$ on the linear span of $\PPP(\MM_1)$. More explicitly, $\Phi$ is defined as follows
    \begin{align*}
        \Phi(\sum_{i=1}^n \lambda_i Q_i): = \sum_{i=1}^n \lambda_i \phi(Q_i),
    \end{align*}
    where $\lambda_i \in \C$ and $Q_i \in \PPP_{1/2}(\MM)$. Let $W := V \otimes E_{12} + V^* \otimes E_{21}$. Note that
    \begin{align*}
        \cos(t)I \otimes (E_{11} -E_{22}) + \sin(t)W
    \end{align*}
 is a self-adjoint unitary such that
 \begin{align*}
     \tau_1(\cos(t)I \otimes (E_{11} -E_{22}) + \sin(t)W) = 0,
 \end{align*}
where $\tau_1$ is the tracial state on $\MM_1$. Therefore $\cos(t)I \otimes (E_{11} -E_{22}) + \sin(t)W$ is in the linear span of $\PPP_{1/2}(\MM_1)$, and
\begin{align*}
    \Phi(\cos(t)I \otimes (E_{11} -E_{22}) + \sin(t)W) = \cos(t)I \otimes (E_{11} -E_{22}) + \sin(t)\Phi(W)
\end{align*}
    is a self-adjoint unitary for every $t \in \R$. In particular $\Phi(W)$ is a self-adjoint unitary. Since $[\cos(t)I \otimes (E_{11} -E_{22}) + \sin(t)\Phi(W)]^2 = I$, we have
\begin{align*}
    [I \otimes (E_{11}-E_{22})] \Phi(W) [I \otimes (E_{11}-E_{22})] = -\Phi(W).
\end{align*}
Therefore their exists a unitary $\psi(V) \in \NN_2$ such that
\begin{align*}
    \Phi(W) = \psi(V) \otimes E_{12} + \psi(V)^* \otimes E_{21}.
\end{align*}
Then it is easy to check that the equation (\ref{equ:map_pos}) holds.
\end{proof}

\begin{lemma}\label{lem:fix_con_fix_all}
    Let $\MM = \NN \otimes M_2(\C)$ where $\NN$ is a finite factor. If $\phi$ is a $L^p$-isometry on $\PPP_{1/2}(\MM)$ such that $\phi(I \otimes E_{11}) = I \otimes E_{11}$ and
    \begin{align*}
        \phi(\frac{1}{2}(I \otimes I_2 +  V \otimes E_{12} + V^* \otimes E_{21})) =
        \frac{1}{2}(I \otimes I_2 +  V \otimes E_{12} + V^* \otimes E_{21})
    \end{align*}
    for very unitary $V \in \NN$, then $\phi(Q) = Q$ for every $Q \in \PPP_{1/2}(\MM)$.
\end{lemma}

\begin{proof}
    There exist a unitary $W_0$ and a positive contraction $H$ in $\NN$ such that
    \begin{align*}
        Q = H \otimes E_{11} + \sqrt{H(I-H)}W_0 \otimes E_{12} + W_0^*\sqrt{H(I-H)} \otimes E_{21} + W_0^*(I-H)W_0 \otimes E_{22}.
    \end{align*}
Let $W := \frac{1}{\sqrt{2}}(I \otimes E_{11} + iW_0 \otimes E_{12} -iI\otimes E_{21} - W_0 \otimes E_{22})$. It is easy to check that $W$ is a unitary and
\begin{align*}
    Q = \frac{1}{2}W^*(I \otimes I_2 + V_0 \otimes E_{12} + V_0^* \otimes E_{21})W,
\end{align*}
where $V_0 = -2\sqrt{H(I-H)} -i(I-2H)$ is a unitary in $\NN$.

Let $\phi_W$ be the surjective $L^p$-isometry on $\PPP_{1/2}(\MM)$ defined as follows
    \begin{align*}
        \phi_W(E) := W\phi(W^*EW)W^*, \quad \forall E \in \PPP_{1/2}(\MM).
    \end{align*}
    It is easy to check that $\phi_W(P) = P$,
    \begin{align*}
     \phi_W(\frac{1}{2}(I \otimes I_2 + iI \otimes E_{12} - iI \otimes E_{21})) = I \otimes I_2 + iI \otimes E_{12} - iI \otimes E_{21},
    \end{align*}
    and
    \begin{align*}
    \phi_W(\frac{1}{2}(I \otimes I_2 + U \otimes E_{12} + U \otimes E_{21})) = I \otimes I_2 + U \otimes E_{12} + U \otimes E_{21}
\end{align*}
    for every self-adjoint unitary $U \in \NN$. By \cref{lem:map_pos} and \cref{thm:unitary}, we know that
\begin{align*}
 \phi_W(\frac{1}{2}(I \otimes I_2 + V \otimes E_{12} + V \otimes E_{21})) = I \otimes I_2 + V \otimes E_{12} + V \otimes E_{21},
\end{align*}
for every unitary $V \in \NN$. In particular,
\begin{align*}
    \phi(Q) = W^*\phi_{W}(\frac{1}{2}(I \otimes I_2 + V_0 \otimes E_{12} + V_0 \otimes E_{21}))W = Q.
\end{align*}
\end{proof}

\begin{theorem}\label{thm:projection}
Let $\MM$ be a finite factor and $\phi$ be a surjective $L^p$-isometry on $\PPP_{1/2}(\MM)$. Then there exists a $*$-automorphism or a *-anti-automorphism $\varphi$ on $\MM$ such that $\phi$ is of one of the following forms:
\begin{enumerate}
    \item $\phi(P) = \varphi(P)$, $P \in \PPP_{1/2}(\MM)$;
    \item $\phi(P) = I- \varphi(P)$, $P \in \PPP_{1/2}(\MM)$.
\end{enumerate}
\end{theorem}

\begin{proof}
    Without loss of generality, we could assume that $\MM = \NN \otimes M_2(\C)$ where $\NN$ is a finite factor and $\phi(I \otimes E_{11}) = I \otimes E_{11}$. By \cref{lem:map_pos}, there exists a surjective $L^p$-isometry $\psi$ on $\UUU(\NN)$ such that
    \begin{align*}
        \phi(\frac{1}{2}(I \otimes I_2 +  V \otimes E_{12} + V^* \otimes E_{21})) =
        \frac{1}{2}(I \otimes I_2 +  \psi(V) \otimes E_{12} + \psi(V)^* \otimes E_{21})
    \end{align*}
for every $V \in \UUU(\NN)$. We can further assume that $\psi(I) = I$.

    By \cref{thm:unitary}, there exists a *-automorphism or a *-anti-automorphism $\varphi$ on $\NN$ such that $\psi(V) = \varphi(V)$ for every unitary $V$ or $\psi(V) = \varphi(V^*)$ for every unitary $V$. Let $\sigma$ be the *-anti-automorphism on $M_2(\C)$ defined by $\sigma(A) = A^t$, where $A^t$ is the transpose of $A$. We define a new surjective $L^p$-isometry $\phi_1$ on $\PPP_{1/2}(\MM)$ as follows.
    \begin{enumerate}
        \item If $\varphi$ is a *-automorphism and $\psi(V) = \varphi(V)$ for every unitary $V$, let
            \begin{align*}
                \phi_1(Q) : = (\varphi^{-1} \otimes \id) \circ \phi(Q).
            \end{align*}
        \item If $\varphi$ is a *-automorphism and $\psi(V) = \varphi(V^*)$ for every unitary $V$, let
            \begin{align*}
                \phi_1(Q) : = I \otimes I_2 - [I \otimes (E_{12}-E_{21})] (\varphi^{-1} \otimes \id)\circ \phi(Q)[I \otimes (E_{21}-E_{12})].
            \end{align*}
        \item If $\varphi$ is a *-anti-automorphism and $\psi(V) = \varphi(V)$ for every unitary $V$, let
            \begin{align*}
                \phi_1(Q) : = I \otimes I_2 - [I \otimes (E_{12}-E_{21})] (\varphi^{-1} \otimes \sigma) \circ \phi(Q) [I \otimes (E_{21}-E_{12})].
            \end{align*}
        \item If $\varphi$ is a *-anti-automorphism and $\psi(V) = \varphi(V^*)$ for every unitary $V$, let
            \begin{align*}
                \phi_1(Q) : = (\varphi^{-1} \otimes \sigma) \circ \phi(Q).
            \end{align*}
    \end{enumerate}
    It is not hard to check that $\phi_1(I \otimes E_{11}) = I \otimes E_{11}$ and
\begin{align*}
    \phi_1(\frac{1}{2}(I \otimes I_2 +  V \otimes E_{12} + V^* \otimes E_{21})) =
        \frac{1}{2}(I \otimes I_2 +  V \otimes E_{12} + V^* \otimes E_{21})
\end{align*}
for every unitary $V \in \NN$. Then \cref{lem:fix_con_fix_all} implies that $\phi_1(Q)=Q$ for every $Q \in \PPP_{1/2}(\MM)$. And the theorem is proved.
\end{proof}

Let $\MM_i$ be a semifinite factor equipped with a n.s.f. traical weight $\tau_i$, $i=1,2$. Assume that there exists a surjective $L^p$-isometry $\phi:\PPP_c(\MM_1) \to \PPP_c(\MM_2)$. Recall that $\phi$ preserves orthogonality in both directions. Therefore, $\MM_1$ is properly infinite if and only if $\MM_2$ is properly infinite. Moreover, it can be shown that $\MM_1$ and $\MM_2$ are either *-isomorphic or *-anti-isomorphic if there exists a surjective $L^p$-isometry $\PPP_c(\MM_1) \to \PPP_c(\MM_2)$.

\begin{lemma}\label{lem:iso_infty_case}
    Let $\MM_1$ and $\MM_2$ be two properly infinite semifinite factors, i.e., $\tau_i(I) = \infty$. If there exists a surjective $L^p$-isometry $\phi:\PPP_c(\MM_1)\rightarrow\PPP_c(\MM_2)$, then there exists a *-isomorphism or a *-anti-isomorphism $\varphi: \MM_1 \to \MM_2$ such that $\tau_1 = \tau_2 \circ \varphi$.
\end{lemma}

\begin{proof}
Let $E_1, E_2 \in\PPP_c(\MM_1)$ such that $E_1\perp E_2$. Since $\phi$ is surjective and preserves orthogonality in both directions, there exists a Hilbert space $\HHH$ such that
   \begin{align*}
      \MM_1 \simeq E_1 \MM_1 E_1 \otimes \BBB(\HHH), \quad \MM_2 \simeq \phi(E_1) \MM_1 \phi(E_1) \otimes \BBB(\HHH).
   \end{align*}
To proved the lemma, we only need to show that there exists a *-isomorphism or a *-anti-isomorphism $\varphi_0: E_1 \MM_1 E_1 \to \phi(E_1) \MM_2 \phi(E_1)$. Invoking \cref{lemma2.6-3-12}, we obtain
    \begin{align*}
        \phi(\{E \in \PPP_c(\MM_1): E \perp (E_1+E_2)\}) = \{E' \in \PPP_c(\MM_2): E' \perp (\phi(E_1)+ \phi(E_2))\}.
    \end{align*}
    Therefore $\phi(F) < \phi(E_1) + \phi(E_2)$ for every projection $F \in \PPP_c(\MM_1)$ satisfying $F < E_1 + E_2$. And $\phi$ restrict to a surjective $L^p$-isometry from $\PPP_{1/2}((E_1+E_2)\MM_1(E_1+E_2))$ to $\PPP_{1/2}((\phi(E_1)+\phi(E_2))\MM_2(\phi(E_1)+\phi(E_2)))$. By \cref{lem:map_pos} and \cref{thm:unitary}, there exists a *-isomorphism or a *-anti-isomorphism $\varphi_0: E_1 \MM_1 E_1 \to \phi(E_1) \MM_2 \phi(E_1)$.
\end{proof}

\begin{proof}[\bf{Proof of \cref{thm:_prop_infty_case}}]
    By \cref{lem:iso_infty_case}, we may assume that $\MM_1 = \MM_2$ and $\tau_1 = \tau_2$. Recall that $\phi$ preserves the orthogonality in both direction. Then the theorem is an immediate consequence of Theorem 1.2 in \cite{Ge1} and Theorem 4.9 in \cite{Shi}.
\end{proof}

\begin{lemma}\label{lem:iso_finite_case}
Let $\MM_i$ be a finite factor with the tracial state $\tau_i$, $i=1,2$. If there exists a surjective $L^p$-isometry $\phi: \PPP_c(\MM_1) \to \PPP_c(\MM_2)$, then $\MM_1$ and $\MM_2$ are either *-isomorphic or *-anti-isomorphic.
\end{lemma}

\begin{proof}
    We only need to prove the lemma under the assumption that $c \in (0, \frac{1}{2}]$. It is well-known that $\MM_1$ and $\MM_2$ are *-(anti)-isomorphic if there exist $s \in (0, 1)$ and nonzero projections $P_1 \in \PPP_s(\MM_1)$ and $P_2 \in \PPP_s(\MM_2)$ such that $P_1 \MM_1 P_1$ and $P_2 \MM_2 P_2$ are *-(anti)-isomorphic. Let $P \in \PPP_c(\MM_1)$. If $c \in (0, \frac{1}{3}] \cup \{\frac{1}{2}\}$, then an argument similar to the one used in the proof of \cref{lem:iso_infty_case} shows that $P\MM_1 P$ and $\phi(P)\MM_2 \phi(P)$ are either *-isomorphic or *-anti-isomorphic.

We now assume that $c \in (\frac{1}{3}, \frac{1}{2})$. Let $f$ be the continuous function on $(-\infty, 1)$ defined by $f(x) = \frac{1-2x}{1-x}$. By \cref{lemma2.6-3-12}, we have
    \begin{align*}
        \phi(\{Q \in \PPP_c(\MM_1): QP = 0\}) = \{Q' \in \PPP_c(\MM_2): Q'\phi(P) = 0\}.
    \end{align*}
Let $c_1 = f(c)$. Note that $0 < c_1 < \frac{1}{2}$. It is clear that $E \in \PPP_{c_1}((I-P)\MM_1 (I-P))$ if and only if $E \in \PPP_{1-2c}(\MM_1)$ and $E < I-P$. And the map $\phi_1$ defined as follows is a surjective $L^p$-isometry from $\PPP_{c_1}((I-P)\MM_1 (I-P))$ to $\PPP_{c_1}((I-\phi(P))\MM_2 (I-\phi(P)))$:
    \begin{align*}
        \phi_1(E) := I-\phi(P) - \phi(I-P - E).
    \end{align*}
If $c_1 \notin (\frac{1}{3}, \frac{1}{2})$, then the discussion above implies that $E \MM_1 E$ and $\phi_1(E) \MM_2 \phi_1(E)$ are either *-isomorphic or *-anti-isomorphic.

Note that the only fixed point of $f$ in $(\frac{1}{3}, \frac{1}{2})$ is $\frac{3-\sqrt{5}}{2}$. We claim that if $c \neq \frac{3-\sqrt{5}}{2}$, then there exists $n \in \N$ such that $c_n:=f^n(c) \notin (\frac{1}{3}, \frac{1}{2})$, where $f^n$ is the n-fold composition of $f$ with itself. We prove the claim by contradiction. Assume that $c \neq \frac{3-\sqrt{5}}{2}$ and $c_n \in (\frac{1}{3}, \frac{1}{2})$ for every $n \in N$. Since
    \begin{align*}
        \begin{cases}
            x - f^2(x) = x + 1/x - 3 > 0, & \forall x \in (\frac{1}{3}, \frac{3-\sqrt{5}}{2})\\
            x - f^2(x) = x + 1/x - 3 < 0, & \forall x \in (\frac{3-\sqrt{5}}{2}, \frac{1}{2})\\
        \end{cases}
    \end{align*}
$\{c_{2k}\}_k$ is a monotonic sequence in $[\frac{1}{3}, \frac{1}{2}]$. Then $\lim_{k} c_{2k}$ is a fixed point of $f^2$. Since $\lim_{k} c_{2k} \neq \frac{3-\sqrt{5}}{2}$ and $\frac{3-\sqrt{5}}{2}$ is the only fixed point of $f^2$ in $[\frac{1}{3}, \frac{1}{2}]$, we have a contradiction and the claim is proved.

If $c \neq \frac{3-\sqrt{5}}{2}$, let $n$ be the number such that $c_n \notin (\frac{1}{3}, \frac{1}{2})$ and $c_i \in (\frac{1}{3}, \frac{1}{2})$ for every $i < n$, then we can repeat the above argument $n$-times and show that there exist $s \in (0,1)$ and $P_1 \in \PPP_{s}(\MM_1)$ and $P_2 \in \PPP_{s}(\MM_2)$ such that $P_1\MM_1 P_2$ and $P_2\MM_2 P_2$ are either *-isomorphic or *-anti-isomorphic.

From now on $c = \frac{3-\sqrt{5}}{2}$. Note that $\MM_1$ and $\MM_2$ are II$_1$ factors since $\frac{3-\sqrt{5}}{2}$ is an irrational number. Let $P \in \PPP_{c}(\MM_1)$. By the same argument used in the proof of Proposition 4.5 in \cite{Shi}, we can show that $(I-P)\MM_1(I-P)$ and $(I-\phi(P))\MM_2(I-\phi(P))$ are either *-isomorphic or *-anti-isomorphic. We therefore only sketch the proof and refer the interested reader to \cite{Shi} for the complete argument.

    Let $\II_{1}=\{E \in\PPP(\MM_1):\tau_1(E)\geq c, E \leq I-P\}$ and $\II_{2}=\{F \in\PPP(\MM_2):\tau_2(F)\geq c, F \leq I-\phi(P)\}$. The map $\phi^+: \II_1 \to \II_2$ defined as follows is a bijection which preserves orthogonality and order in both directions (see Proposition 3.1 in \cite{Shi} for the proof):
\begin{align*}
    \phi^+(E):=\vee\{\phi(Q):Q\in\PPP_c(\MM_1), Q\leq E\}.
\end{align*}
Note that $\phi^+(I-P) = I-P$. We have
   \begin{align*}
       \phi'(E):=I-\phi(P)-\phi^+(I-P-E)
   \end{align*}
    is a map from $\{E \in \PPP(\MM_1): \tau_1(E) \leq 1-2c, E < I-P\}$ to $\{F \in \PPP(\MM_2): \tau_2(F) \leq 1-2c, F < I-\phi(P)\}$. Proceed as in the proof of Proposition 4.5 in \cite{Shi}, it can be shown that $\phi'$ is a bijection which preserves orthogonality and order in both directions. By remark 3.3 in \cite{Shi}, $\phi'$ can be extended to an ortho-isomorphism from $\PPP((I-P)\MM_1(I-P))$ to $\PPP((I-\phi(P))\MM_2(I-\phi(P)))$. Then Corollary in \cite{Dye} implies that $(I-P)\MM_1(I-P)$ and $(I-\phi(P))\MM_2(I-\phi(P))$ are either *-isomorphic or *-anti-isomorphic.
\end{proof}

\begin{proof}[\bf{Proof of \cref{thm:finite_case}}]
    We only need to prove the theorem under the further assumption that $c \in (0, 1/2]$. By \cref{lem:iso_finite_case}, we may assume that $\MM_1 = \MM_2$ and $\tau_1 = \tau_2$. Recall that $\phi$ preserves the orthogonality in both direction. Then the theorem is an immediate consequence of Theorem 1.2 in \cite{Ge1}, Theorem 4.9 in \cite{Shi} and \cref{thm:projection}.
\end{proof}


\begin{thebibliography}{99}
\bibitem{AE} E. Andruchow, Short geodesics of unitaries in the $L^2$-metric, Can. Math. Bull., 48(3) (2005), 340-354.

\bibitem{BJM}
    F. Botelho, J. Jamison and L. Moln\'{a}r,
    Surjective isometries on Grassmann spaces,
    J. Funct. Anal. 265(2013), 2226-2238.

\bibitem{XBC} T. N. Bekjan, Z. Chen and Q. Xu, Introduction to operator algebras and noncommutative $L^p$-spaces (in Chinese), (2010) Science Press (China).

\bibitem{Dye} H. A. Dye, On the geometry of projections in certain operator algebras, Ann. Math., 61(1955), 73-89.

\bibitem{Ge1} G. P. Geh$\acute{\mathrm{e}}$r, P. $\breve{\mathrm{S}}$emrl, Isometries of Grassmann spaces, J. Funct. Anal., 270 (2016), 1585-1601.

 \bibitem{Ge2} G. P. Geh$\acute{\mathrm{e}}$r, Wigner's theorem on Grassmann spaces, J. Funct. Anal., 273 (2017), 2994-3001.

\bibitem{Ge3} G. P. Geh$\acute{\mathrm{e}}$r, P. $\breve{\mathrm{S}}$emrl, Isometries of Grassmann spaces, II, Adv. Math., 332 (2018), 287-310.

\bibitem{Shi} W. Gu, M. Ma, J. Shen and W. Shi, Ortho-isomorphisms on Grassmann spaces in factors of type II, to appear in Journal of Operator Theory.

%\bibitem{Gy} M. Gy\"{o}ry, Transformations on the set of all n-dimensional subspaces of a Hilbert space preserving orthogonality, Publ. Math. Debrecen 65(2004), 233-242

\bibitem{HLP} G. H. Hardy, J. E. Littlewood and G. P\'{o}lya, Inequalities. Cambridge University Press, Cambridge (1934).

\bibitem{HHMM} O. Hatori, G. Hirasawa, T. Miura and L. Moln$\acute{\mathrm{a}}$r, Isometries and maps compatible with inverted Jordan triple products on groups, Tokyo J. Math., (35) 2 (2012), 385-410.

\bibitem{HM} O. Hatori, L. Moln$\acute{\mathrm{a}}$r, Isometries of the unitary groups, Proc. Amer. Math. Soc., 130 (2012), 2141-2154.

\bibitem {Lam} J. Lamperti, On the isometries of certain function spaces, Pacific Journal of Mathematics, (8) 3 (1958), 459-466.

\bibitem{Kad} R. Kadison, Isometries of operator algebras, Ann. Math., (2) 54 (1951), 325-338.

 \bibitem{M1} L. Moln$\acute{\mathrm{a}}$r, Generalization of Wigner’s unitary-antiunitary theorem for indefinite inner product spaces, Commun. Math. Phys., 201(2000), 785-791.

 \bibitem{M2} L. Moln$\acute{\mathrm{a}}$r, Transformations on the set of all n-dimensional subspaces of a Hilbert space preserving principal angles, Commun. Math. Phys., 217(2001), 409-421.


\bibitem {LM} L. Moln$\acute{\mathrm{a}}$r, Jordan triple endomorphisms and isometries of unitary groups,  Linear Algebra Appl., (439) 11 (2013), 3518-3531.

\bibitem{MS} L. Moln$\acute{\mathrm{a}}$r, P. $\breve{\mathrm{S}}$emrl, Transformations of the unitary group on a Hilbert space, J. Math. Anal. Appl., 388(2012), 1205-1217.

 \bibitem{MM} M. Mori, Isometries between projection lattices of von Neumann algebras, J. Funct. Anal., 276(2019), 3511-3528.

%\bibitem{Pa} M. Pankov, Geometric version of Wigner's theorem for Hilbert Grassmannians, J. Math. Anal. Appl., 459  (2018) 135-144.

\bibitem{Pan} M. Pankov, Wigner-typetheorems for Hilbert Grassmannians,LMS Lecture Note Series 460,Cambridge University Press, Cambridge (2020).

\bibitem{QWWY} W. Qian, L. Wang, W. Wu and W. Yuan, Wigner-type theorem on transition probability preserving maps in semfinite factors, J. Funct. Anal., 276(2019), 1773-1787.

\bibitem {Se} P. $\breve{\mathrm{S}}$emrl, Orthogonality preserving transformations on the set of $n$-dimensional subspaces of a Hilbert space, Illinois J. Math., 48 (2004), 567-573.

%\bibitem{Uh} U. Uhlhorn, Representation of symmetry transformations in quantum mechanics, Ark. Fysik 23(1963), 307-340.

\bibitem{Tak}
    M. Takesaki, Theory of Operator Algebras I, II, III, Encyclopaedia of Mathematical Sciences, Springer-Verlag Berlin Heidelberg (2002, 2003).

\bibitem{Wi}E. P. Wigner, Gruppentheorie und ihre Anwendung auf die Quanten mechanik der Atomspektren, Fredrk Vieweg und Sohn, 1931(English translation Group Theory and its Applications to the Quantum Mechanics of Atomic Spectra, Academic Press, 1959).
\end{thebibliography}
\end{document}